\documentclass[11pt,a4paper,reqno]{amsart}
\usepackage{amssymb,amsmath,amsthm}
\usepackage[utf8]{inputenc}
\usepackage[T1]{fontenc}
\usepackage{enumerate}
\usepackage[all]{xy}
\usepackage{fullpage}
\usepackage{comment}
\usepackage{array}
\usepackage{longtable}
\usepackage{stmaryrd}
\usepackage{mathrsfs} 
\usepackage{xcolor}
\usepackage{mathtools}

\def\gcd{\mathrm{gcd}}
\def\deg{\mathrm{deg}}

\def\dim{\mathrm{dim}}

\usepackage{hyperref}
\hypersetup{
    colorlinks=true,
    linkcolor=blue,
    filecolor=red,  
    citecolor=green,
    urlcolor=cyan,
    pdftitle={GON},
    pdfpagemode=FullScreen,
    }

\usepackage{cleveref}
\crefformat{section}{\S#2#1#3}
\crefformat{subsection}{\S#2#1#3}
\crefformat{subsubsection}{\S#2#1#3}
\usepackage{enumitem}
\usepackage{tikz}
\usepackage{mathdots}

\newtheorem{theorem}{Theorem}[section]
\newtheorem{lemma}[theorem]{Lemma}

\newtheorem{corollary}[theorem]{Corollary}
\newtheorem{proposition}[theorem]{Proposition}
\newtheorem{example}[theorem]{Example}

\theoremstyle{definition}
\newtheorem{definition}[theorem]{Definition}

\theoremstyle{remark}
\newtheorem{remark}[theorem]{Remark}

\makeatletter
\newcommand{\subalign}[1]{%
  \vcenter{%
    \Let@ \restore@math@cr \default@tag
    \baselineskip\fontdimen10 \scriptfont\tw@
    \advance\baselineskip\fontdimen12 \scriptfont\tw@
    \lineskip\thr@@\fontdimen8 \scriptfont\thr@@
    \lineskiplimit\lineskip
    \ialign{\hfil$\m@th\scriptstyle##$&$\m@th\scriptstyle{}##$\hfil\crcr
      #1\crcr
    }%
  }%
}
\makeatother

\numberwithin{equation}{section}

\title{On the Minimal Denominator Problem in Function Fields}
\author{Noy Soffer Aranov}
\email{noysofferaranov@math.utah.edu}
\address{Department of Mathematics, University of Utah, Salt Lake City, Utah, USA}

\begin{document}

\maketitle
\begin{abstract}
    We study the minimal denominator problem in function fields. In particular, we compute the probability distribution function of the the random variable which returns the degree of the smallest denominator $Q$, for which the ball of a fixed radius around a point contains a rational function of the form $\frac{P}{Q}$. Moreover, we discuss the distribution of the random variable which returns the denominator of minimal degree, as well as higher dimensional and $P$-adic generalizations. This can be viewed as a function field generalization of a paper by Chen and Haynes. 
\end{abstract}
\paragraph{MSC Class: 11J04, 11J13, 11J61}

\section{Introduction}
Meiss and Sanders \cite{MS} described an experiment in which a distance $\delta>0$ is fixed, and for randomly chosen $x\in [0,1)$, they study the statistics of the function 
\begin{equation}
    q_{\min}(x,\delta)=\min\left\{q:\exists\frac{p}{q}\in B(x,\delta),\gcd(p,q)=1\right\}.
\end{equation}
Chen and Haynes \cite{CH} computed the the probability that $\mathbb{P}(q_{\min}(x,\delta)=q)$ for every $\delta>0$ and for every $q\leq \left[\frac{1}{\delta}\right]$. Moreover, they proved that $\mathbb{E}[q_{\min}(\cdot, \delta)]=\frac{16}{\pi^2\cdot \delta^{\frac{1}{2}}}+O(\log^2\delta)$. Artiles \cite{Art} and Markloff \cite{M} generalized the results of \cite{CH} to higher dimensions by studying the statistics of Farey fractions and using dynamical methods. The minimal denominator problem was investigated in the real setting in several other papers such as \cite{KM,St}, but it is not well studied over other fields. 

In this paper, we use linear algebra and number theory to study the function field analogue of the function $q_{\min}(x,\delta)$, as well as its higher dimensional and $P$-adic analogues in the function field setting. In particular, we prove a function field analogue of the results of \cite{CH}. We note that unlike \cite{CH,M,Art}, we do not study the distribution of Farey fractions or use dynamical techniques, rather we use linear algebra and lattice point counting techniques, which provide better bounds in ultrametric spaces. 
\subsection{Function Field Setting}
In this setting, let $q$ be a prime power, and denote the ring of Laurent polynomials over $\mathbb{F}_q$ by 
$$\mathcal{R}=\left\{\sum_{n=0}^Na_nx^n:a_n\in \mathbb{F}_q,N\in \mathbb{N}\cup\{0\}\right\}.$$
Let $\mathcal{K}$ be the field of fractions of $\mathcal{R}$, and define an absolute value on $\mathcal{K}$ by $\left|\frac{f}{g}\right|=q^{\deg(f)-\deg(g)}$, where $f,g\in \mathcal{R}$ and $g\neq 0$. Then, the completion of $\mathcal{K}$ with respect to $\vert \cdot\vert$ is
$$\mathcal{K}_{\infty}=\left\{\sum_{n=-N}^{\infty}a_nx^{-n}:a_n\in \mathbb{F}_q\right\}.$$
Let $\mathcal{O}=\{\alpha\in \mathcal{K}_{\infty}:\vert \alpha\vert\leq 1\}$, and let 
$$\mathfrak{m}=x^{-1}\mathcal{O}=\{\alpha\in \mathcal{K}_{\infty}:\vert \alpha\vert\leq q^{-1}\}.$$
For $\alpha\in \mathcal{K}_{\infty}$, we write $\alpha=[\alpha]+\{\alpha\}$, where $[\alpha]\in \mathcal{R}$ and $\{\alpha\}\in \mathfrak{m}$. In this paper, the Haar measure on $\mathcal{K}_{\infty}$ is defined to be the unique translation invariant measure $\mu$, such that $\mu(\mathfrak{m})=1$. 

In $\mathcal{K}_{\infty}^n$, the supremum norm is defined as $\Vert (v_1,\dots,v_n)\Vert=\max_{i=1,\dots,n}\Vert \mathbf{v}_i\Vert$. Similarly, for $\boldsymbol{\alpha}=(\alpha_1,\dots,\alpha_n)\in \mathcal{K}_{\infty}^n$, we let $[\boldsymbol{\alpha}]=([\alpha_1],\dots,[\alpha_n])$ and $\{\boldsymbol{\alpha}\}=(\{\alpha_1\},\dots,\{\alpha_n\})$. 
\subsection{Main Results}
We prove a function field analogue of the main results of \cite{CH}. Let $n\in \mathbb{N}$. For $\delta>0$ and $\boldsymbol{\alpha}\in\mathcal{K}_{\infty}^n$, we define the minimal denominator degree by
$$\deg_{\min}(\boldsymbol{\alpha},\delta)=\min\left\{d:\exists\frac{P}{Q},\deg(Q)=d,\left|\boldsymbol{\alpha}-\frac{P}{Q}\right|<\delta\right\}.$$
We say that $Q$ is a minimal denominator for $\boldsymbol{\alpha}$ if $\deg(Q)=\deg_{\min}(\boldsymbol{\alpha},\delta)$ and $\left\Vert\boldsymbol{\alpha}-\frac{P}{Q}\right\Vert<\delta$. We note that if $Q$ is a minimal denominator for $\boldsymbol{\alpha}$, then, it is also a minimal denominator for $\{\boldsymbol{\alpha}\}$. Hence, it suffices to assume that $\boldsymbol{\alpha}\in \mathfrak{m}^n$. Moreover, since the absolute value $\vert \cdot \vert$ obtains values in $\{0\}\cup\{q^{k}:k\in \mathbb{Z}\}$, then, for every $q^{-(k+1)}<\delta\leq q^{-k}$, we have $\deg_{\min}(\boldsymbol{\alpha},\delta)=\deg_{\min}(\boldsymbol{\alpha},q^{-k})$. Hence, it suffices to take $\delta=q^{-k}$, where $k\in \mathbb{N}$. We first compute the probability distribution function of $\deg_{\min}(\cdot,q^{-k})$ when $n=1$. From now on, we denote the probability distribution by $\mathbb{P}$.
\begin{theorem}
\label{thm:deg_min1D}
    Let $k\in \mathbb{N}$. Then, we have
    $$\mathbb{P}\left(\deg_{\min}(\alpha,q^{-1})=d\right)=\begin{cases}
        \frac{1}{q}& \text{ if } d=0,\\
        \frac{q-1}{q}&\text{ if }d=1
    \end{cases},$$
    and for every $k\geq 2$, we have
    \begin{equation}
        \mathbb{P}\left(\deg_{\min}(\alpha,q^{-k})=d\right)=\begin{cases}
            q^{-k}&\text{ if }d=0,\\
            \frac{q-1}{q^{k-2d+1}}&\text{ if }d\leq \left\lceil\frac{k}{2}\right\rceil,d\in \mathbb{N},\\
            0&\text{ else}.
        \end{cases}
    \end{equation}
\end{theorem}
\begin{corollary}
    Let $k\in \mathbb{N}$. Then,
    \begin{equation}
        \mathbb{E}[\deg_{\min}(\cdot,q^{-k})]=\begin{cases}
            \frac{q-1}{q}&\textit{ if }k=1,\\
            \frac{q-1}{q^{k-1}}\left(\frac{q^{2\left\lceil\frac{k}{2}\right\rceil+1}\left(\left\lceil\frac{k}{2}\right\rceil+1\right)-\left(\left\lceil\frac{k}{2}\right\rceil+2\right)q^{2\left\lceil\frac{k}{2}\right\rceil}+1}{(q^2-1)^2}\right)&\text{else}.
        \end{cases}
    \end{equation}
\end{corollary}
\begin{proof}
    When $k=1$, the claim follows immediately from Theorem \ref{thm:deg_min1D}. Otherwise, by Theorem \ref{thm:deg_min1D}, we have
    \begin{equation}
    \begin{split}
        \mathbb{E}\left[\deg_{\min}(\alpha,q^{-k})\right]=\sum_{d=0}^{\left\lceil\frac{k}{2} \right\rceil}d\frac{q-1}{q^k}q^{2d-1}=\frac{q-1}{q^{k-1}}\frac{d}{dt}\left(\sum_{d=0}^{\left\lceil\frac{k}{2}\right\rceil}t^d\right)_{t=q^2}\\
        =\frac{q-1}{q^{k-1}}\frac{d}{dt}\left(\frac{t^{\left\lceil\frac{k}{2}\right\rceil+1}-1}{t-1}\right)_{t=q^2}=\frac{q-1}{q^{k-1}}\left(\frac{q^{2\left\lceil\frac{k}{2}\right\rceil+1}\left(\left\lceil\frac{k}{2}\right\rceil+1\right)-\left(\left\lceil\frac{k}{2}\right\rceil+2\right)q^{2\left\lceil\frac{k}{2}\right\rceil}+1}{(q^2-1)^2}\right).
    \end{split}
    \end{equation}
\end{proof}
Moreover, in every dimension, there is a unique monic polynomial which is a denominator of minimal degree. 
\begin{lemma}
\label{lem:UniqueQ_min}
    For every $\boldsymbol{\alpha}\in \mathfrak{m}^n$ and for every $k\geq 1$, there exists a unique monic polynomial $Q\in \mathcal{R}$, such that $\deg(Q)=\deg_{\min}(\boldsymbol{\alpha},q^{-k})$ and $\Vert Q\boldsymbol{\alpha}\Vert<q^{-k}$.
\end{lemma}
This motivates the following definition.
\begin{definition}
    Due to Lemma \ref{lem:UniqueQ_min}, we denote the unique monic polynomial $Q$ satisfying $\deg(Q)=\deg_{\min}(\alpha,q^{-k})$ and $\Vert Q\alpha\Vert<q^{-k}$ by $Q_{\min}(\alpha,q^{-k})$. 
\end{definition}
We also compute the distribution of $Q_{\min}(\cdot,q^{-k})$. To do so, we shall use some notations from number theory.
\begin{definition}
    For a polynomial $Q$, we let $d(Q)$ be the number of prime divisors of $Q$, we let $D(Q)$ be the number of monic divisors of $Q$, and we let $S(Q)$ be the set of divisors of $Q$. Define
    $$\mu(Q)=\begin{cases}
        (-1)^{d(Q)}&\text{if }Q\text{ is square free},\\
        0&\text{if there exists }P\text{ such that }P^2\mid Q
    \end{cases}$$
\end{definition}
\begin{definition}
For a polynomial $Q\in \mathcal{R}$, we define $S_{\text{monic}}^{*,\ell}(Q)$ to be the set of $\ell$ tuples $(a_1,\dots,a_{\ell})$, such that $a_i$ are distinct monic polynomials which divide $Q$, and $\deg(a_i)<\deg(Q)$. 
\end{definition}
\begin{theorem}
\label{thm:Q_min=Q}
    Let $Q$ be a monic polynomial with $\deg(Q)\leq \left\lceil\frac{k}{2}\right\rceil$. Then, for every $k\geq 1$, the probability that $Q_{\min}(\alpha,q^{-k})=Q$ is
    \begin{equation}
    \begin{split}
        \frac{1}{q^k}\left(\vert Q\vert+\sum_{N|Q,\deg(N)<\deg(Q)}\vert N\vert\sum_{\ell=1}^{D(N)}(-1)^{\ell}\left(\frac{D\left(\frac{Q}{N}\right)!}{\left(D\left(\frac{Q}{N}\right)-\ell\right)!}+\sum_{M\in S\left(\frac{Q}{N}\right):D\left(\frac{Q}{NM}\right)\geq \ell}\mu(M)\frac{D(M)!}{(D(M)-\ell)!}\right)\right).
    \end{split}
    \end{equation}
    In particular, if $Q$ is an irreducible monic polynomial of degree $d$, then, 
    \begin{equation}
        \mathbb{P}\left(Q_{\min}(\alpha,q^{-k})=Q\right)=\frac{q^d-1}{q^k}.
    \end{equation}
\end{theorem}
The higher dimensional setting discusses a simultaneous solution for the equations $\vert \{Q\alpha_i\}\vert<q^{-k}$.
\begin{lemma}
\label{lem:deg_min>=maxdeg}
    For every $k\in \mathbb{N}$, and for every $\boldsymbol{\alpha}=(\alpha_1,\dots,\alpha_n)\in \mathfrak{m}^{n}$, we have 
    $$\deg_{\min}(\boldsymbol{\alpha},q^{-k})\geq \max_{i=1,\dots,n} \deg_{\min}(\alpha_i,q^{-k}).$$
\end{lemma}
\begin{proof}
    If there exists $Q\in \mathcal{R}$ with $\deg(Q)=\deg_{\min}(\boldsymbol{\alpha},q^{-k})$, such that for every $i=1,\dots,n$, we have $\Vert Q\alpha_i\Vert<q^{-k}$, then, $\deg(Q)\geq \deg_{\min}(\alpha_i,q^{-k})$ for every $i=1,\dots,n$. Hence, $\deg_{\min}(\boldsymbol{\alpha},q^{-k})\geq \max \deg_{\min}(\alpha_i,q^{-k})$.
\end{proof}
Hence, it is natural to ask what is the probability that $\deg_{\min}(\boldsymbol{\alpha},q^{-k})=\max_{i=1,\dots,n}\deg_{\min}(\alpha_i,q^{-k})$. This is an immediate corollary of Lemma \ref{lem:UniqueQ_min} and Lemma \ref{lem:deg_min>=maxdeg}
\begin{corollary}
\label{cor:Q_min=maxQ}
    Let $\boldsymbol{\alpha}\in \mathfrak{m}^n$ and let $k\geq 1$. Then, for $Q\in \mathcal{R}$ with $\deg(Q)=d$, with $0\leq d\leq \left\lceil\frac{k}{2}\right\rceil$, we have $Q_{\min}(\boldsymbol{\alpha},q^{-k})=Q$ if and only if 
    \begin{enumerate}
        \item \label{cond:Q_min|Q_i}for every $i=1,\dots, n$, we have $Q_{\min}(\alpha_i,q^{-k})|Q$;
        \item \label{cond:Q_min,i=Q}there exists $i$ such that $Q_{\min}(\alpha_i,q^{-k})=Q$. 
    \end{enumerate}
\end{corollary}
\begin{proof}[Proof of Corollary \ref{cor:Q_min=maxQ}]
    Assume that \eqref{cond:Q_min|Q_i} and \eqref{cond:Q_min,i=Q} hold. For $i=1,\dots,n$, let $\frac{P_i}{Q_i}\in \mathcal{K}$ be such that $Q_{\min}(\alpha_i,q^{-k})=Q_i$ and $\left|\alpha_i-\frac{P_i}{Q_i}\right|<q^{-k}$. Since $Q_i|Q$, then, there exists $D_i\in \mathcal{R}$, such that $Q=Q_iD_i$. Therefore, $\left|\alpha_i-\frac{P_iD_i}{Q}\right|=\left|\alpha_i-\frac{P_i}{Q_i}\right|<q^{-k}$, so that $\left\Vert\boldsymbol{\alpha}-\frac{\mathbf{P}}{Q}\right\Vert<q^{-k}$, where $\mathbf{P}=(P_1,\dots,P_n)$. Thus, $\deg_{\min}(\boldsymbol{\alpha},q^{-k})\leq \deg(Q)$. On the other hand, by Lemma \ref{lem:deg_min>=maxdeg}, $\deg_{\min}(\boldsymbol{\alpha},q^{-k})=\deg(Q)$. Hence, by Lemma \ref{lem:UniqueQ_min} and \eqref{cond:Q_min,i=Q}, $Q_{\min}(\boldsymbol{\alpha},q^{-k})=Q$.
    
    Assume that $Q_{\min}(\boldsymbol{\alpha},q^{-k})=Q$ and that there exists $\mathbf{P}=(P_1,\dots,P_n)$ such that $\left\Vert \boldsymbol{\alpha}-\frac{\mathbf{P}}{Q}\right\Vert<q^{-k}$. By Lemma \ref{lem:deg_min>=maxdeg}, we have $\deg_{\min}(\boldsymbol{\alpha},q^{-k})\geq \deg_{\min}(\alpha_i,q^{-k})$ for every $i=1,\dots,n$. For $i=1,\dots,n$, let $Q_i=Q_{\min}(\alpha_i,q^{-k})$ and let $P_i'$ be such that $\left|\alpha_i-\frac{P_i'}{Q_i}\right|<q^{-k}$. If $Q_i=Q$ for every $i$, there is nothing to prove. Otherwise, let $i$ be such that $Q_i\neq Q$. By Lemma \ref{lem:UniqueQ_min}, $\left|\alpha_i-\frac{P_i'}{Q_i}\right|<q^{-k}$ and $\left|\alpha_i-\frac{P_i}{Q}\right|<q^{-k}$ imply together that $\frac{P_i}{Q}=\frac{P_i'}{Q_i}$. Thus, $Q_i$ divides $Q$.
\end{proof}
Furthermore, we bound the probability that $\deg_{\min}(\boldsymbol{\alpha},q^{-k})\leq d$ in higher dimensions. 
\begin{lemma}
\label{lem:multdimUppBnd}
    Let $n,k\in \mathbb{N}$. Then, for every $\boldsymbol{\alpha}\in \mathfrak{m}^n$, we have $\deg_{\min}(\boldsymbol{\alpha},q^{-k})\leq \left\lceil\frac{nk}{n+1}\right\rceil$. Moreover, for every $d\leq \left\lceil\frac{nk}{n+1}\right\rceil$, we have
    \begin{equation}
        \mathbb{P}(\deg_{\min}(\boldsymbol{\alpha},q^{-k})\leq d)\leq q^{-(kn-(n+1)d)}.
    \end{equation}
\end{lemma}
We also discuss a $P$-adic variant of the minimal denominator problem, which is motivated by the $P$-adic Littlewood conjecture \cite{dMT,EK,BKN,ALN,Rob,GR}. Let $P\in \mathcal{R}$ be an irreducible polynomial, $\alpha\in \mathfrak{m}$, and let $k\geq 1$. Define
$$\deg_{\min,P}(\alpha,q^{-k})=\min\left\{d\geq 0:\exists m\geq 0, \exists\frac{a}{b}:\deg(b)=d, \frac{a}{P^mb}\in B(\alpha,q^{-k})\right\}=\inf_{m\geq 0}\deg_{\min}(P^m\alpha,q^{-k}).$$
\begin{theorem}
\label{thm:Padic}
    For every irreducible polynomial $P$ and for almost every $\alpha\in \mathfrak{m}$, we have $\deg_{\min,P}(\alpha,q^{-k})=0$. 
\end{theorem}
Theorem \ref{thm:Padic} motivates us to study the set of $\alpha\in \mathfrak{m}$ such that $\deg_{\min,P}(\alpha,q^{-k})\geq 1$ using subtler notions of "size". One such notion to study the exceptional set is Hausdorff dimension. We first recall the definition of Hausdorff dimension in the function field setting. 
\begin{definition}
    Let $E\subseteq \mathcal{K}_{\infty}^n$, let $\delta>0$, and let $s\in \mathbb{R}$. Define 
    $$H_{\delta}^{(s)}(E)=\inf\left\{\sum_{i=1}^{\infty}r_i^s:E\subseteq \bigcup_{i\in \mathbb{N}}B(\boldsymbol{\alpha}_i,r_i),r_i\leq \delta\right\},$$
    where $B(\boldsymbol{\alpha},r)=\{\boldsymbol{\beta}\in \mathcal{K}_{\infty}^n:\Vert \boldsymbol{\alpha}-\boldsymbol{\beta}\Vert\leq r\}$. For $s\in \mathbb{R}$, let $H^{(s)}(E)=\lim_{\delta\rightarrow 0}H_{\delta}^{(s)}(E)$, and define 
    $$\dim_H(E)=\inf\{s:H^{(s)}(E)=0\}.$$
\end{definition}
Hausdorff dimension can often be difficult to compute, but for the purposes of this paper, it suffices to understand the following example. Similar examples, which are often called Cantor sets or missing digit sets, appear in books in fractal geometry such as \cite{Fal}.
\begin{example}
\label{ex:MissingWord}
    Let $P\in \mathbb{F}_q[x]$ be a polynomial and let $D\subseteq (\mathbb{F}_q[x]/P\mathbb{F}_q[x])^k$. Let $E$ be the set of $\alpha=\sum_{n=-\infty}^Na_nP^n\in \mathcal{K}_{\infty}$, where $N\in \mathbb{N}$ and $a_n\in \mathbb{F}_q[x]/P\mathbb{F}_q[x]$, such that for every $n\in \mathbb{Z}$, we have $(a_1,\dots,a_n)\in D$. Then, $\dim_H(E)=\frac{\log_q\vert D\vert}{\log_q\vert P\vert^k}=\frac{\log_q\vert D\vert}{k\deg(P)}$.
\end{example}
\begin{theorem}
    \label{thm:PadicDim}
    For every irreducible $P\in \mathcal{R}$ and for every $1\leq d\leq \left\lceil\frac{k}{2}\right\rceil$, we have
    $$\dim_H\{\alpha\in \mathfrak{m}:\deg_{\min,P}(\alpha,q^{-k})\geq d\}=\frac{\log_q(\vert P\vert^k-\vert P\vert^{2(d-1)})}{k\deg(P)}.$$
\end{theorem}
\begin{remark}
A natural question pertains to the minimal denominator question in the infinite residue case, for example over $\mathbb{Q}[x]$ or $\mathbb{C}[x]$. Since many of our results rely on counting techniques, these results do not hold for function fields with an infinite residue field.
\end{remark}
\subsection{Acknowledgements}
I would like to thank Eran Igra and Albert Artiles who introduced me to the minimal denominator problem, and thus led to the birth of this paper. I would also like to thank the anonymous referee, whose comments helped improve this paper. 
\section{Hankel Matrices}
\label{sec:HankelMatrix}
We first translate the minimal denominator problem to the language of linear algebra. For $k,\ell\in \mathbb{N}$ and $\alpha=\sum_{i=1}^{\infty}\alpha_ix^{-i}\in \mathfrak{m}$, we define the Hankel matrix of $\alpha$ of order $(k,\ell)$ as
$$\Delta_{\alpha}(k,\ell)=\begin{pmatrix}
    \alpha_1&\alpha_2&\dots&\alpha_{\ell}\\
    \alpha_2&\alpha_3&\dots&\alpha_{\ell+1}\\
    \vdots&\dots&\ddots&\vdots\\
    \alpha_k&\alpha_{k+1}&\dots&\alpha_{k+\ell-1}
\end{pmatrix}.$$
Assume that $\alpha\in \mathfrak{m}$ and $\frac{P}{Q}\in \mathbb{F}_q(x)$. Then, $\left|\alpha-\frac{P}{Q}\right|<q^{-k}$ if and only if $\vert Q\alpha-P\vert<q^{-(k-d)}$, where $d=\deg(Q)$. Using Dirichlet's Theorem in function fields (see \cite{dMT70}, \cite[Theorem 1.1]{GaGh17} and \cite[Theorem 4.1]{A}), one can trivially bound $\deg_{\min}(\alpha,q^{-k})$.
\begin{theorem}{Dirichlet's Theorem in Function Fields \cite[Theorem 1.1]{GaGh17}}
\label{thm:Dirichlet}
    For every $\alpha\in\mathfrak{m}$ and every $k\in \mathbb{N}$, there exist $P\in \mathcal{R}$ and $Q\in \mathcal{R}\setminus\{0\}$, such that 
    \begin{equation}
    \begin{cases}
        \left|Q\alpha-P\right|<q^{-k}, \textit{and}\\
        \vert Q\vert\leq q^{k}.
    \end{cases}
    \end{equation}
\end{theorem}
Note that Theorem \ref{thm:Dirichlet} implies that $\deg_{\min}(\alpha, q^{-k})\leq k$. Hence, we can assume that $d\leq k$. Then, $\vert\{ Q\alpha\}\vert<q^{-(k-d)}$ if and only if
\begin{equation}
\label{eqn:HankelMinDenom}
    \begin{pmatrix}
        \alpha_1&\alpha_2&\dots&\alpha_{d+1}\\
        \alpha_2&\alpha_3&\dots&\alpha_{d+2}\\
        \vdots&\dots&\ddots&\vdots\\
        \alpha_{k-d}&\alpha_{k-d+1}&\dots&\alpha_k
    \end{pmatrix}\begin{pmatrix}
        Q_0\\
        Q_1\\
        \vdots\\
        Q_d
    \end{pmatrix}=0,
\end{equation}
where $Q=\sum_{i=0}^dQ_ix^i$. We notice that if $d+1\geq k-d+1$, that is if $d\geq \frac{k}{2}$, then, there exists a non-trivial solution to (\ref{eqn:HankelMinDenom}). Hence, $\deg_{\min}(\alpha, q^{-k})\leq \left\lceil\frac{k}{2}\right\rceil$. We first use this reinterpretation to prove Lemma \ref{lem:UniqueQ_min}. To do so we define for $\boldsymbol{\alpha}\in \mathfrak{m}^n$ the matrix
    $$\Delta_{\boldsymbol{\alpha}}(k,\ell)=\begin{pmatrix}
        \Delta_{\alpha_1}(k,\ell)\\
        \vdots\\
        \Delta_{\alpha_n}(k,\ell)
    \end{pmatrix}.$$
\begin{lemma}
\label{rem:degRank}
    We have $\deg_{\min}(\boldsymbol{\alpha},q^{-k})=d$, for $d\leq \left\lceil\frac{k}{2}\right\rceil$, if and only if for every $j<d$, the matrix $\Delta_{\boldsymbol{\alpha}}(k-j,j+1)$ has rank $j+1$, but the matrix $\Delta_{\boldsymbol{\alpha}}(k-d,d+1)$ has rank $d$. In particular $\dim\left(\operatorname{Ker}\left(\Delta_{\boldsymbol{\alpha}}(k-d,d+1)\right)\right)=1$.
\end{lemma}
\begin{proof}
    Notice that if $d=\deg_{\min}(\boldsymbol{\alpha},q^{-k})$, then, there exists $Q=Q_0+Q_1x+\dots+Q_dx^n$ with $Q_j\in \mathbb{F}_q$ and $Q_d\neq 0$, such that $\vert \{Q\alpha_i\}\vert<q^{-(k-d)}$ for every $i=1,\dots,n$. Therefore, 
    \begin{equation}
    \label{eqn:Delta(k-d,d+1)}
        \Delta_{\boldsymbol{\alpha}}(k-d,d+1)\begin{pmatrix}
            Q_0\\
            Q_1\\
            \vdots\\
            Q_d
        \end{pmatrix}=0.
    \end{equation}
    Hence, $\operatorname{rank}(\Delta_{\boldsymbol{\alpha}}(k-d,d+1))\leq d$. On the other hand, since $d=\deg_{\min}(\boldsymbol{\alpha},q^{-k})$, for every non-zero polynomial $P=P_0+P_1x+\dots+P_{j}x^{j}$ with $j\leq d-1$, we have $\Vert \{P\boldsymbol{\alpha}\}\Vert\geq q^{-k}$. Thus, there exists $i=1,\dots,n$, such that $\vert \{P\alpha_i\}\vert\geq q^{-k}$, so that $\Delta_{\alpha_i}(k-j,j+1)(P_0,\dots,P_{d-1})^T\neq 0$ for every non-zero $P$. As a consequence, for every non-zero $P$, we have
    \begin{equation}
    \label{eqn:Delta(k-d-1,d)}
        \Delta_{\boldsymbol{\alpha}}(k-j,j+1)\begin{pmatrix}
            P_0\\
            P_1\\
            \vdots\\
            P_{j}
        \end{pmatrix}\neq 0.
    \end{equation}
    Hence, $\operatorname{rank}(\Delta_{\boldsymbol{\alpha}}(k-j,j+1))=d$. Thus, by \eqref{eqn:Delta(k-d-1,d)}, for every $j=1,\dots,d-1$, we have 
    $$\operatorname{rank}(\Delta_{\boldsymbol{\alpha}}(k-j,j+1))=j+1$$ 
    Furthermore, \eqref{eqn:(k-d,d+1)Sum} and \eqref{eqn:Delta(k-d-1,d)} imply that $d=\operatorname{rank}(\Delta_{\boldsymbol{\alpha}}(k-d,d+1))$ and $\dim\left(\operatorname{Ker}\left(\Delta_{\boldsymbol{\alpha}}(k-d,d+1)\right)\right)=1$.
\end{proof}
We first prove that for every $\alpha$, there is a unique minimal denominator. 
\begin{lemma}
    For every $0\neq \alpha\in \mathfrak{m}$, there exists a unique $d\leq \frac{k}{2}$ for which there exist coprime $P,Q$ with $\deg(Q)=d$ such that $\left|\alpha-\frac{P}{Q}\right|<q^{-k}$.
\end{lemma}
\begin{proof}
    Assume that there exist $P,Q,P',Q'\in \mathcal{R}$ such that 
    \begin{enumerate}
        \item $Q,Q'\neq 0$
        \item $\gcd(P,Q)=1=\gcd(P',Q')$, 
        \item $\deg(Q')=d'<d=\deg(Q)\leq \frac{k}{2}$, 
        \item $\left|\alpha-\frac{P}{Q}\right|<q^{-k}$, and
        \item $\left|\alpha-\frac{P'}{Q'}\right|<q^{-k}$. 
    \end{enumerate}Then, by the ultrametric inequality and the fact that these fractions are reduced, we have
    \begin{equation}
        \left|\frac{PQ'-P'Q}{QQ'}\right|= \left|\frac{P}{Q}-\frac{P'}{Q'}\right|<\frac{1}{q^k}.
    \end{equation}
    Thus, $\vert PQ'-P'Q\vert<\frac{q^{d+d'}}{q^k}\leq q^{2d}{q^{-k}}=q^{-(k-2d)}\leq 1$. Hence, $\frac{P}{Q}=\frac{P'}{Q'}$, which contradicts the assumption that $\operatorname{gcd}(P',Q')=\operatorname{gcd}(P,Q)=1$.
\end{proof}
\begin{proof}[Proof of Lemma \ref{lem:UniqueQ_min}]
    Assume that $\deg_{\min}(\boldsymbol{\alpha},q^{-k})=d$. Then, by Lemma \ref{rem:degRank}, $\dim\operatorname{Ker}(\Delta_{\boldsymbol{\alpha}}(k-d,d+1))=1$. Let $Q$ be a monic polynomial satisfying $\deg(Q)=d$ and $\Vert Q\boldsymbol{\alpha}\Vert<q^{-k}$. Thus, $$\Delta_{\boldsymbol{\alpha}}(k-d,d+1)(Q_0,Q_1,\dots,Q_{d-1},1)^T=0.$$
    Since $\dim\operatorname{Ker}(\Delta_{\boldsymbol{\alpha}}(k-d,d+1))=1$, then, $(Q_0,Q_1,\dots,Q_{d-1},1)$ is the unique vector $\mathbf{v}=(v_0,\dots,v_d)$ with $v_d=1$, such that $\Delta_{\boldsymbol{\alpha}}(k-d,d+1)\mathbf{v}=0$. Thus, by (\ref{eqn:HankelMinDenom}), $Q$ is the unique monic polynomial of minimal degree with $\Vert Q\boldsymbol{\alpha}\Vert<q^{-k}$.
\end{proof}
We shall use several facts about ranks of Hankel matrices to prove Theorems \ref{thm:deg_min1D} and \ref{thm:Q_min=Q}. 
\begin{theorem}{\cite[Theorem 5.1]{AGR}}
\label{thm:numHankMatrix}
    Let $r>0$. Then, the number of invertible $h\times h$ Hankel matrices with entries in $\mathbb{F}_q$ of rank $r$, $N(r,h;q)$, is equal to
    \begin{equation}
        N(r,h;q)=\begin{cases}
            1,&\textit{if }r=0,\\
            q^{2r-2}(q^2-1),&\textit{if }1\leq r\leq h-1, \textit{ and}\\
            q^{2h-2}(q-1),&\textit{if }r=h
        \end{cases}.
    \end{equation}
\end{theorem}
We shall also use the following generalization of Theorem \ref{thm:numHankMatrix}.
\begin{theorem}{\cite[Theorem 1.1]{DG}}
    \label{thm:DG}
    Let $k,\ell\in \mathbb{N}$, let $F$ be a finite field with $\vert F\vert=q$, and let $r\leq \min\{k,\ell\}-1$. Then, the number of Hankel matrices $\Delta_{\alpha}(k,\ell)$ over $F$ with rank at most $r$ is $q^{2r}$.
\end{theorem}
\begin{lemma}{\cite[Lemma 2.3]{ALN}}
\label{lem:ALN}
    Let $m,n\in \mathbb{N}$, and let $k\leq \min\{m,n-1\}$. Let $H=\Delta_{\alpha}(m,n)$ be a Hankel matrix. If the first $k$ columns of $H$ are independent, but the first $k+1$ columns of $H$ are dependent, then, $\det(\Delta_{\alpha}(k,k))\neq 0$.
\end{lemma}
\section{Proofs in the One Dimensional Case} 
\begin{proof}[Proof of Theorem \ref{thm:deg_min1D}]
By using the reinterpretation in \cref{sec:HankelMatrix}, we realize that if $k=1$, then, \eqref{eqn:HankelMinDenom} has a non-trivial solution when $d\geq 1$. Moreover, there exists a solution for (\ref{eqn:HankelMinDenom}) when $k=1$ and $d=0$ if and only if $\alpha_1=0$. Hence, 
$$\mathbb{P}\left(\deg_{\min}(\alpha,q^{-1})=0\right)=\frac{1}{q}, \mathbb{P}\left(\deg_{\min}(\alpha,q^{-1})=1\right)=\frac{q-1}{q}.$$
If, $k\geq 2$ and $d\geq \frac{k}{2}$, then, there exists a non-trivial solution to (\ref{eqn:HankelMinDenom}). Hence, $\deg_{\min}(\alpha,q^{-k})\leq \left\lceil\frac{k}{2}\right\rceil$. Firstly, if $d=0$, then, $\alpha_1=\dots=\alpha_k$, and therefore, $\mathbb{P}(\deg_{\min}(\alpha,q^{-k})=0)=q^{-k}$.

Let $1 \leq d\leq \frac{k}{2}$ and let $\alpha\in \mathbb{F}_q((x^{-1}))$. By Lemma \ref{rem:degRank}, we have $d=\deg_{\min}(\alpha,q^{-k})$ if and only if the columns of the Hankel matrix $\Delta_{\alpha}(k-d,d+1)$ are linearly dependent, but the columns of $\Delta_{\alpha}(k-d-1,d)$ are linearly independent. Hence, by Lemma \ref{lem:ALN}, the matrix $\Delta_{\alpha}(d,d)$ is invertible. Hence, there exist unique $a_1,\dots, a_d\in \mathbb{F}_q$, such that
\begin{equation}
\label{eqn:(d,d)MinorSum}
    \begin{pmatrix}
        \alpha_{d+1}\\
        \alpha_{d+2}\\
        \vdots\\
        \alpha_{2d}
    \end{pmatrix}=a_1\begin{pmatrix}
        \alpha_1\\
        \alpha_2\\
        \vdots\\
        \alpha_d
    \end{pmatrix}+\dots+a_d\begin{pmatrix}
        \alpha_d\\
        \alpha_{d+1}\\
        \vdots\\
        \alpha_{2d-1}
    \end{pmatrix}.
\end{equation}
On the other hand, since the columns of $\Delta_{\alpha}(k-d,d+1)$ are linearly dependent, and the columns of $\Delta_{\alpha}(k-d+1,d)$ are linearly independent, there exist $b_1,\dots,b_d\in \mathbb{F}_q$, such that
\begin{equation}
\label{eqn:(k-d,d+1)Sum}
    \begin{pmatrix}
        \alpha_{d+1}\\
        \alpha_{d+2}\\
        \vdots\\
        \alpha_k
    \end{pmatrix}=b_1\begin{pmatrix}
        \alpha_1\\
        \alpha_2\\
        \vdots\\
        \alpha_{k-d}
    \end{pmatrix}+\dots+b_d\begin{pmatrix}
        \alpha_d\\
        \alpha_{d+1}\\
        \vdots\\
        \alpha_{k-1}
    \end{pmatrix}.
\end{equation}
Thus, by (\ref{eqn:(d,d)MinorSum}) and (\ref{eqn:(k-d,d+1)Sum}), we have $a_i=b_i$ for every $i=1,\dots, d$. Hence, given an invertible matrix $d\times d$ Hankel matrix $\Delta_{\alpha}(d,d)$ and some $\alpha_{2d}\in \mathbb{F}_q$, there is exactly one way to extend the word $(\alpha_1,\dots,\alpha_{2d})$ to a word of length $k$, $(\alpha_1,\dots,\alpha_k)$, such that for every choice of $\{\alpha_j\}_{j\geq k+1}\subseteq \mathbb{F}_q$, the Laurent sequence $\sigma=\sum_{i=1}^{\infty}\alpha_ix^{-i}$ satisfies $\deg_{\min}(\sigma,q^{-k})=d$. Therefore, by Theorem \ref{thm:numHankMatrix} (see also Theorem \ref{thm:DG}), we have
$$\mathbb{P}(\deg_{\min}(\alpha,q^{-k})=d)=\frac{q^{2d-1}(q-1)}{q^k}.$$

\end{proof}
We now prove Theorem \ref{thm:Q_min=Q}. In order to do so, we utilize the following counting claim.
\begin{proposition}
\label{prop:GCDCnt}
    The number of primitive vectors in $S_{\text{monic}}^{*,\ell}(Q)$ is
    \begin{equation}
        \vert\widehat{S}_{\text{monic}}^{*,\ell}\vert=\begin{cases}
            \frac{D(Q)!}{(D(Q)-\ell)!}+\sum_{N\in S(Q):D\left(\frac{Q}{N}\right)\geq \ell}\mu(N)\frac{D(N)!}{(D(N)-\ell)!} & D(Q)\geq \ell,\\
            0& \text{else}.
        \end{cases}
    \end{equation}
\end{proposition}
\begin{remark}
    We use Proposition \ref{prop:GCDCnt} to count the number of tuples whose greatest common denominator is $1$, instead of the classical counting method, since this yields a more compact expression in the proof of Theorem \ref{thm:Q_min=Q}
\end{remark}
\begin{proof}[Proof of Proposition \ref{prop:GCDCnt}]
    We first note that 
    $$\vert S_{\text{monic}}^{*,\ell}(Q)\vert=\begin{cases}\frac{D(Q)!}{(D(Q)-\ell)!}&\vert D(Q)\vert\geq \ell,\\
    0&\text{else}.
    \end{cases}$$
    Hence, to count primitive vectors in $S_{\text{monic}}^{*,\ell}$, we use the inclusion-exclusion principle and induction to obtain that
    \begin{equation}
    \begin{split}
        \vert \widehat{S}_{\text{monic}}^{*,\ell}\vert=\left|S_{\text{monic}}^{*,\ell}\setminus\bigcup_{P\in S(Q)\text{ prime}}PS_{\text{monic}}^{*,\ell}\left(\frac{Q}{P}\right)\right|\\
        =\vert S_{\text{monic}}^{*,\ell}\vert-\sum_{P_1,\dots,P_i\in S(Q)\text{ prime}}(-1)^{i+1}\left|P_1\cdots P_iS_{\text{monic}}^{*,\ell}\left(\frac{Q}{P_1\cdots P_i}\right)\right|\\
        =\frac{D(Q)!}{(D(Q)-\ell)!}+\sum_{N\in S(Q)}\mu(N)\left|S_{\text{monic}}^{*,\ell}\left(\frac{Q}{N}\right)\right|\\
        =\frac{D(Q)!}{(D(Q)-\ell)!}+\sum_{N\in S(Q):D\left(\frac{Q}{N}\right)\geq \ell}\mu(N)\frac{D(N)!}{(D(N)-\ell)!}.
    \end{split}
    \end{equation}
\end{proof}
\begin{proof}[Proof of Theorem \ref{thm:Q_min=Q}]
Let $Q$ be a monic polynomial of degree at most $\left\lceil\frac{k}{2}\right\rceil$. We use the following fact from \cite{CR}: We have $\Vert Q\alpha\Vert<q^{-(k-d)}$, where $Q=Q_0+Q_1x+\dots+Q_{d-1}x^{d-1}+x^d$, if and only if
\begin{equation}
\label{eqn:CircMatrix}
    \begin{pmatrix}
        Q_0&\dots&Q_{d-1}&-1&0&\dots&0\\
        0&Q_0&\dots&Q_{d-1}&-1&0&\dots\\
        \vdots&\dots&\ddots&\ddots&\dots&\ddots&\ddots\\
        0&\dots&\dots&Q_0&\dots&Q_{d-1}&-1
    \end{pmatrix}\begin{pmatrix}
        \alpha_1\\
        \alpha_2\\
        \vdots\\
        \alpha_{k}
    \end{pmatrix}=\boldsymbol{0}.
\end{equation}
We denote the matrix the left hand side of (\ref{eqn:CircMatrix}) by $A_Q$, and we denote $\pi_k(\alpha)=(\alpha_1,\dots,\alpha_k)^t$. Then, $A_Q\in M_{k-d\times k}(\mathbb{F}_q)$, and $\dim(\operatorname{Ker}(A_Q))=d$. Hence, $\vert \operatorname{Ker}(A_Q)\vert=q^d$. By Lemma \ref{lem:UniqueQ_min}, if $\pi_k(\alpha)\in \operatorname{Ker}(A_D)\cap \operatorname{Ker}(A_Q)$, where $\deg(D)<\deg(Q)$, then, $D|Q$. Hence, by the inclusion-exclusion principle,
\begin{equation}
\begin{split}
\label{eqn:P(Q_min=Q)}
\mathbb{P}(Q_{\min}(\alpha,q^{-k})=Q)=\mathbb{P}\left(\pi_k(\alpha)\in \left(\operatorname{Ker}(A_Q)\right)\setminus\bigcup_{D|Q,\deg(D)<\deg(Q)}\operatorname{Ker}(A_D)\right)\\
=\mathbb{P}(\pi_k(\alpha)\in \operatorname{Ker}(A_Q))-\sum_{D_1,\dots,D_{\ell}|Q}(-1)^{\ell+1}\mathbb{P}\left(\bigcap_{i=1}^{\ell}\operatorname{Ker}(A_{D_i})\right)\\
=q^{-(k-\deg(Q))}+\sum_{D_1,\dots,D_{\ell}|Q}(-1)^{\ell}\mathbb{P}\left(\operatorname{Ker}(A_{\gcd(D_1,\dots,D_{\ell}})\right).
\end{split}
\end{equation}
We notice that if $N|Q$, then, $N=\gcd(D_1,\dots,D_{\ell})$ if and only if for every $i=1,\dots,\ell$, there exists a monic polynomial $a_i\in \mathcal{R}$ such that $D_i=a_iN$, $a_i\mid \frac{Q}{N}$, and $\gcd(a_1,\dots,a_{\ell})=1$. Hence, $(a_1,\dots,a_{\ell})\in \mathcal{R}^{\ell}$ is a primitive vector with distinct coordinates, which are all monic polynomials dividing $\frac{Q}{N}$, so that $(a_1,\dots,a_{\ell})\in S_{\text{monic}}^{*,\ell}\left(\frac{Q}{N}\right)$. Hence, by Proposition \ref{prop:GCDCnt}, we have
\begin{equation}
\begin{split}
    \sum_{D_1,\dots,D_{\ell}|Q}(-1)^{\ell+1}\mathbb{P}\left(\operatorname{Ker}(A_{\operatorname{gcd}(D_1,\dots,D_{\ell})})\right)\\
    =\frac{1}{q^k}\sum_{N|Q,\deg(N)<\deg(Q)}\vert N\vert\sum_{\ell=1}^{D(Q)}(-1)^{\ell+1}\#\{(D_1,\dots,D_{\ell}):\gcd(D_1,\dots,D_{\ell})=N\}
    \\=\frac{1}{q^k}\sum_{N|Q,\deg(N)<\deg(Q)}\vert N\vert\sum_{\ell=1}^{D(N)}(-1)^{\ell+1}\left|\widehat{S}_{\text{monic}}^{*,\ell}\left(\frac{Q}{N}\right)\right|\\
    =\frac{1}{q^k}\sum_{N|Q,\deg(N)<\deg(Q)}\vert N\vert\sum_{\ell=1}^{D(N)}(-1)^{\ell+1}\left(\frac{D\left(\frac{Q}{N}\right)!}{\left(D\left(\frac{Q}{N}\right)-\ell\right)!}+\sum_{M\in S\left(\frac{Q}{N}\right):D\left(\frac{Q}{N}\right)\geq \ell}\mu(M)\frac{D(M)!}{(D(M)-\ell)!}\right).
\end{split}
\end{equation}
Hence, the probability that $Q_{\min}(\alpha,q^{-k})=Q$, for $\deg(Q)\leq \left\lceil\frac{k}{2}\right\rceil$ is equal to
\begin{equation*}
\begin{split}
    \frac{1}{q^k}\left(\vert Q\vert+\sum_{N|Q,\deg(N)<\deg(Q)}\vert N\vert\sum_{\ell=1}^{D(N)}(-1)^{\ell}\left(\frac{D\left(\frac{Q}{N}\right)!}{\left(D\left(\frac{Q}{N}\right)-\ell\right)!}+\sum_{M\in S\left(\frac{Q}{N}\right):D\left(\frac{Q}{NM}\right)\geq \ell}\mu(M)\frac{D(M)!}{(D(M)-\ell)!}\right)\right).
\end{split}
\end{equation*}
\end{proof}
\section{Proofs in the $P$-adic and High Dimensional Cases}
\begin{proof}[Proof of Theorem \ref{thm:Padic}]
We notice that $\left|P^mb\alpha-a\right|<q^{-(k-d)}$, where $d=\deg(b)$, if and only if 
\begin{equation}
    \begin{pmatrix}
        \alpha_m&\alpha_{m+1}&\dots&\alpha_{m+d}\\
        \alpha_{m+1}&\alpha_{m+2}&\dots&\alpha_{m+d+1}\\
        \vdots&\dots&\ddots&\vdots\\
        \alpha_{m+k-d}&\dots&\dots&\alpha_{m+k}
    \end{pmatrix}\begin{pmatrix}
        b_0\\
        b_1\\
        \vdots\\
        b_d
    \end{pmatrix}=0,
\end{equation}
where $\alpha=\sum_{i=1}^{\infty}P^{-i}\alpha_i$, for some $\alpha_i\in \mathbb{F}_q[x]/P\mathbb{F}_q[x]$. Hence, $\deg_{\min,P}(\alpha,q^{-k})\leq \frac{k}{2}$. Since almost every string is normal \cite{Bor}, then, for every $k$, for every prime $P$, and for almost every $\alpha\in \mathfrak{m}$, the string $0^k$ appears in the infinite word $\{\alpha_i\}_{i\in \mathbb{N}}$. Hence, if $\alpha_m=\alpha_{m+1}=\dots=\alpha_{m+k-1}=0$, then, there exists $a$, such that $\vert P^m\alpha-a\vert<q^{-(k-d)}$. Thus, $\mathbb{P}(\deg_{\min,P}(\alpha,q^{-k})=0)=1$. 
\end{proof}
\begin{proof}[Proof of Theorem \ref{thm:PadicDim}]
By Theorem \ref{thm:DG}, for every $k$ and every $d\leq k-d$, the number of Hankel matrices $\operatorname{rank}(\Delta_{\alpha}(k-d-1,d))=d$ with entries in $\mathbb{F}_q[x]/P\mathbb{F}_q[x]$ of rank $d$ is $\vert P\vert^{k}-\vert P\vert ^{2(d-1)}$. Thus, by Example \ref{ex:MissingWord}, for every $1\leq d\leq \frac{k}{2}$, we have
$$\dim_H\{\alpha\in \mathfrak{m}:\deg_{\min,P}(\alpha,q^{-k})\geq d\}=\frac{\log_q(\vert P\vert ^k-\vert P\vert^{2(d-1)})}{k\deg(P)}.$$
\end{proof} 
\begin{proof}[Proof of Lemma \ref{lem:multdimUppBnd}]
We note that $\Vert Q\boldsymbol{\alpha}\Vert<q^{-k}$ if and only if 
\begin{equation}\label{eqn:highDimMinDenom}\Delta_{\boldsymbol{\alpha}}(k-d,d+1)(Q_0,Q_1,\dots,Q_d)^T=0.\end{equation} 
Thus, if $d\geq \frac{nk}{n+1}$, (\ref{eqn:highDimMinDenom}) has a non-trivial solution. If $d=0$, then, for every $i=1,\dots,n$, we have $(\alpha_{i1},\dots,\alpha_{ik})=0$, where $\alpha_i=\sum_{j=1}^{\infty}\alpha_{ij}x^{-j}$. Hence, 
\begin{equation}\label{eqncase:0}\mathbb{P}\left(\deg_{\min}(\boldsymbol{\alpha},q^{-k})=0\right)=q^{-nk}.\end{equation}
Let $1\leq d\leq \left\lceil\frac{nk}{n+1}\right\rceil$, and let $Q$ be a monic polynomial with $\deg(Q)=d$. Since $\vert \operatorname{Ker}(A_Q)\vert=q^d$, we have
\begin{equation}
    \mathbb{P}\left(Q_{\min}(\boldsymbol{\alpha},q^{-k})=Q\right)\leq \mathbb{P}\left(\bigcap_{i=1}^n\left(\alpha_i\in \operatorname{Ker}(A_Q)\right)\right)=\prod_{i=1}^n\mathbb{P}\left(\alpha_i\in \operatorname{Ker}(A_Q)\right)=q^{-n(k-d)}
\end{equation}
Hence, 
\begin{equation}
    \mathbb{P}\left(\deg_{\min}(\boldsymbol{\alpha},q^{-k})=d\right)\leq q^{-nk+(n+1)d}.
\end{equation}
\end{proof} 
\bibliography{Ref}
\bibliographystyle{amsalpha}
\end{document}